\documentclass[letterpaper, 10 pt, conference]{ieeeconf}  

\IEEEoverridecommandlockouts                              
\usepackage{amsmath, amssymb, amsxtra, mathrsfs}
\usepackage{graphics, epsfig, here, fancyhdr,graphicx}
\usepackage{subfigure}

\newtheorem{theorem}{Theorem}
\newtheorem{definition}[theorem]{Definition}

\newcommand{\bR}{\mathbb{R}}

\newcommand{\cH}{\mathcal{H}}

\newcommand{\cO}{\mathcal{O}}

\newcommand{\bq}[1]{{\left[#1\right]}}

\newcommand{\norm}[1]{{\left\|#1\right\|}}

\newcommand{\rM}[1]{{\mathrm{#1}}}

\newcommand{\w}{\omega}

\newcommand{\eP}{\varepsilon}

\newcommand{\rd}{\mathrm{d}}
\newcommand{\ddx}[1]{{\frac{\rd }{\rd #1}}}

\title{\LARGE \bf
Synthesis of Optimal Ensemble Controls for Linear Systems using the Singular Value Decomposition$^*$
}

\author{Anatoly Zlotnik$^\dagger$ and Jr-Shin Li$^\ddagger$
\thanks{$^\dagger$azlotnik@ese.wustl.edu, \,\, Electrical and Systems Engineering, Washington University in St. Louis, One Brookings Drive, Saint Louis, MO 63130
        }%
\thanks{$^\ddagger$jsli@ese.wustl.edu, \,\,Electrical and Systems Engineering, Washington University in St. Louis, One Brookings Drive, Saint Louis, MO 63130
        }%
\thanks{$^*$This work was supported by the NSF Career Award \#0747877 and the AFOSR YIP FA9550-10-1-0146.}
}

\begin{document}

\maketitle
\thispagestyle{empty}
\pagestyle{empty}

\begin{abstract}

An emerging and challenging area in mathematical control theory called Ensemble Control encompasses a class of problems that involves the guidance of an uncountably infinite collection of structurally identical dynamical systems, which are indexed by a parameter set, by applying the same open-loop control.  The subject originates from the study of complex spin dynamics in Nuclear Magnetic Resonance (NMR) spectroscopy and imaging (MRI).  A fundamental question concerns ensemble controllability, which determines the existence of controls that transfer the system between desired initial and target states.  For ensembles of finite-dimensional time-varying linear systems, the necessary and sufficient controllability conditions and analytical optimal control laws have been shown to depend on the singular system of the operator characterizing the system dynamics. Because analytical solutions are available only in the simplest cases, there is a need to develop numerical methods for synthesizing these controls. We introduce a direct, accurate, and computationally efficient algorithm based on the singular value decomposition (SVD) that approximates ensemble controls of minimum norm for such systems.  This method enables the application of ensemble control to engineering problems involving complex, time-varying, and high-dimensional linear dynamic systems.

\end{abstract}

\section{Introduction} \label{sec_intro}

The implementation of all scientific and engineering applications is complicated by uncertainty or variation in system model parameters, for which known control techniques are unable to successfully compensate.  This issue is especially challenging when the control task must be accomplished without feedback, whether the control function must transfer a single control system between states of interest without sensitivity to an uncertain parameter set, or steer a possibly uncountable collection of structurally identical systems with variation in common parameters between states that may depend on the parameters.  Investigation of the latter category is motivated by factors that arise in practical applications of quantum control theory to the fields of Nuclear Magnetic Resonance (NMR) spectroscopy and imaging (MRI), and has given rise to a new area of mathematical control theory called ensemble control \cite{li06thesis}.  Rapidly progressing technologies based on quantum theory require the manipulation of very large ensembles of quantum systems on the order of Avogodro's number ($6\times 10^{23}$), whose states cannot be measured during the transfer, and whose dynamics are subject to dispersion in parameters such as frequency.  The performance of the necessary controls must be insensitive to parameter variation across the ensemble, as well as to inhomogeneity in the applied radiofrequency (RF) control field \cite{li09tac,li06pra}.  A long standing problem of significance to NMR requires the design of RF excitations that steer a given quantum ensemble between initial and target states, and whose performance is insensitive to parameter variation \cite{kobzar04,kobzar05,pauly91}.  An acceptable control function must concurrently drive a collection of systems, with identical dynamics but parameter values unknown up to a given range, between desired initial and target states.

Although first motivated by the necessity to control large collections of similar systems, the mathematical devices produced by investigating ensemble control can also be used to approach any open-loop control application in which the system response must be immune to uncertainty in model parameters.  For instance, harmonic oscillators are widely used to approximate periodic phenomena in a variety of scientific and engineering applications where the frequency of oscillation may not be known exactly, but rather is confined to a given range \cite{bartlett02,brumer92}.  Harmonic oscillators often appear in quantum-electrodynamics, and steering such quantum systems using electromagnetic fields is a subject of widespread interest \cite{mirrahimi04}.

The theoretical investigation of ensemble control begins with the notion of ensemble controllability, which determines the existence of controls that achieve various types of state transfers for a system of interest.  It has been shown that a bilinear system evolving on $\rM{SO}(3)$ called the Bloch equations, which models the evolution over time of a sample of nuclear spins, is ensemble controllable \cite{li06pra,li06thesis}.  In addition, the necessary and sufficient conditions for ensemble controllability of finite-dimensional time-varying linear systems for transfers between states in Hilbert space have been recently derived \cite{li09ccc}.  These conditions depend on the singular system of the operator characterizing the system dynamics, which can in turn be used to represent the minimum norm control that accomplishes the transfer as an infinite sum of weighted eigenfunctions.  This method was used to synthesize optimal ensemble controls for a harmonic oscillator system, for which the resulting eigenfunctions are the well-known family of prolate spheroidal wave functions \cite{li11tac}.  This special structure facilitates synthesis of the controls in this special case, as well as the computation of optimal controls with bounded amplitude by solving a constrained convex optimization problem \cite{li11tac}.

The controllability conditions for general, possibly nonlinear, ensemble control problems are presently unknown, and generalized analytical control design methods remain a challenging problem, although analytical solutions exist for a few specific systems.  When the singular system of the relevant integral operator is unavailable in analytical form, or in the case of a nonlinear ensemble system, an alternative is to use a pseudospectral numerical method that translates an optimal control problem in function space into a discrete nonlinear programming problem \cite{li09jcp,ruths10,li11pnas,ross03,gong06}.
This method has been effective for solving a variety of ensemble control problems, but may be difficult to implement for large-scale systems with variation in many parameters because of the computational complexity required for approximating control functions with sufficient accuracy \cite{ruths11}.  Therefore, a need exists for direct and computationally efficient numerical methods for synthesizing ensemble controls that can accomplish various state transfers for a variety of systems.

In this paper, we introduce an accurate, stable, and computationally efficient numerical method based on the singular value decomposition for constructing minimum norm ensemble controls for finite-dimensional time-varying linear systems.  In Section \ref{sec_theo}, we briefly review the theoretical work that forms the foundation of our approach. In Section \ref{sec_method}, we describe our method for numerically approximating the singular system of the Fredholm integral operator of the first kind that characterizes the dynamics of a linear ensemble system, as well as the synthesis of the unique minimum norm control that accomplishes a desired transfer in function space.  In Section \ref{sec_ex}, we revisit the control of harmonic oscillators to demonstrate the effectiveness of our new method for accomplishing complex state transfers, and also examine ensemble control of linear systems with higher dimension and with variation in several parameters.  Finally in Section \ref{sec_conc}, we conclude by discussing the advantages of the method presented here, as well as our future work on developing fast, iterative methods for ensemble control of nonlinear systems.  Such techniques will accelerate the rapidly expanding scope of ensemble control theory and spectral methods by contributing powerful new tools for solving cutting-edge problems in diverse fields from neuroscience to quantum physics.

\section{Ensemble control of linear systems} \label{sec_theo}

The aim of ensemble control is to simultaneously manipulate a continuum of dynamical systems, which are governed by internal and external dynamics that depend on a parameter varying over a compact indexing set, by applying the same open-loop control input to each.  In this section, we review the basic definitions and the fundamental theoretical results that enable ensemble control synthesis for finite-dimensional time-varying linear systems.

Consider a parameterized family of dynamical systems indexed by a parameter $\beta$ varying over a compact set $K$, given by
\begin{align}\label{ode}
\dot{X}(t,\beta)&=A(t,\beta)X(t,\beta)+B(t,\beta)u(t),\\
&  X\in M \subset\mathbb{R}^n, \quad \beta\in K\subset\mathbb{R^d}, \quad u\in U\subset \mathbb{R}^m, \nonumber
\end{align}
where $A(t,\beta)\in\mathbb{R}^{n\times n}$ and $B(t,\beta)\in\mathbb{R}^{n\times m}$ have elements that are real $L_{\infty}$ and $L_2$ functions, respectively, defined on a compact set $D=[0,T]\times K$, and are denoted $A\in L_\infty^{n\times n}(D)$ and $B\in L_2^{n\times m}(D)$.  The ensemble controllability conditions for the system (\ref{ode}) depend on the existence of an open-loop control $u:[0,T]\to U$ that can steer the instantaneous state of the ensemble $X(t,\cdot):K\to M$ between any points of interest in the Hilbert space of functions on $K$.  Let $\mathcal{H}_T=L_2^m[0,T]$ denote the set of $m$-tuples, whose elements are complex vector-valued square-integrable measurable functions defined on $0\leq t\leq T$, with an inner product defined by
\begin{align} \label{inpt}
\langle g,h \rangle_T = \int_0^T g^\dagger(t) h(t) dt,
\end{align}
where $\dagger$ denotes the conjugate transpose.  Similarly, let $\mathcal{H}_K=L_2^n(K)$ be equipped with an inner product
\begin{align} \label{inpb}
\langle p,q \rangle_K = \int_K p^\dagger(\beta) q(\beta) d\mu(\beta),
\end{align}
where $\mu$ is the Lebesgue measure.  With well-defined addition and scalar multiplication, $\mathcal{H}_T$ and $\mathcal{H}_K$ are separable Hilbert spaces, where $||\cdot||_T$ and $||\cdot||_K$ denote their respective induced norms.
\begin{definition} (\textit{Ensemble controllability} \cite{li09ccc})
We say that the family (\ref{ode}) is ensemble controllable on the function space $\mathcal{H}_K$ if for all $\eP>0$, and all $X_0,\, X_F\in \mathcal{H}_K$, there exists $T>0$ and an open loop piecewise-continuous control $u \in \mathcal{H}_T$, such that starting from $X(0,\beta)=X_0(\beta)$, the final state $X(T,\beta)=X_T(\beta)$ satisfies $||X_T-X_F||_K<\eP$.
\end{definition}
In other words, the system (\ref{ode}) is ensemble controllable if it is possible to guide it from $X_0$ to $X_F$ in the space $\cH_K$, where the acceptable range of $T\in(0,\infty)$ may depend on $\eP$, $K$, and $U$.  Necessary and sufficient conditions have been determined for the ensemble controllability of finite-dimensional time-varying linear systems, and are based on the Fredholm integral operator that characterizes the system dynamics \cite{li11tac}.  Given the initial state $X(0,\beta)=X_0(\beta)$ of the system (\ref{ode}), the variation of parameters formula gives rise to the solution
\begin{align} \label{soldet}
X(T,\beta) & =\Phi(T,0,\beta)X_0(\beta) \nonumber \\
& \quad +\int_{0}^{T}\Phi(T,\sigma,\beta)B(\sigma,\beta)u(\sigma)d\sigma,
\end{align}
where $\Phi(T,0,\beta)$ is the transition matrix for the system $\dot{X}(t,\beta)=A(t,\beta)X(t,\beta)$.  Our goal is for the terminal state to equal the target state in the function space $\mathcal{H}_K$, so setting $X(T,\beta)=X_F(\beta)$, pre-multiplying by $\Phi(0,T,\beta)$ and rearranging results in the integral operator equation
\begin{align} \label{critdet}
(Lu)(\beta)=\int_0^{T}\Phi(0,\sigma,\beta)B(\sigma,\beta)u(\sigma)d\sigma=\xi(\beta),
\end{align}
where $\xi(\beta)=\Phi(0,T,\beta)X_F(\beta)-X_0(\beta)$.  The theory of ensemble controllability and the derivation of minimum norm controls can be reduced to the solvability of the above integral equation.  A spectral decomposition, called the singular system, of the operator $L$ is used to produce an infinite eigenfunction series expansion for the $u\in \mathcal{H}_T$ of minimum norm that satisfies (\ref{critdet}) with sufficient accuracy.
\begin{definition} \label{singsys} {\it Singular System} \cite{gohberg03}: Let $Y$ and $Z$ be Hilbert spaces and $L:Y\to Z$ be a compact operator.  If $(\sigma_n^2,\nu_n)$ is an eigensystem of $LL^*$ and $(\sigma_n^2,\mu_n)$ is an eigensystem of $L^*L$, namely, $LL^*\nu_n=\sigma_n^2\nu_n$, $\nu_n\in Z$, and $L^*L\mu_n=\sigma_n^2\mu_n$, $\mu_n\in Y$, where $\sigma_n>0$ ($n\geq 1$), and the two systems are related by the equations $L\mu_n=\sigma_n\nu_n$ and $L^*\nu_n=\sigma_n\mu_n$, we say that $(\sigma_n,\mu_n,\nu_n)$ is a singular system of $L$.
\end{definition}
Suppose that $(\sigma_n,\mu_n,\nu_n)$ is a singular system of the operator $L$ as defined in (\ref{critdet}), which is compact \cite{li09ccc}.  The necessary and sufficient conditions for ensemble controllability of the system (\ref{ode}) have been proven to be
\begin{align} \label{contcond}
 (i) & \quad\sum_{n=1}^{\infty}\frac{|\langle\xi,\nu_n\rangle_T|^2}{\sigma_n^2}<\infty, \\
(ii) & \quad \xi\in\overline{\mathcal{R}(L)},
\end{align}
where $\overline{\mathcal{R}(L)}$ denotes the closure of the range space of $L$.  In addition, it was shown that the control law
\begin{align} \label{cont1}
u=\sum_{n=1}^\infty \frac{\langle \xi, \nu_n\rangle_T}{\sigma_n}\mu_n
\end{align}
satisfies $\langle u,u\rangle_T\leq \langle u_0,u_0\rangle_T$ for all $u_0\in \mathcal{U}$ and $u_0\neq u$, where $\mathcal{U}=\{v\mid Lv=\xi \}$.  For further details we refer the reader to the original work on this subject  \cite{li09ccc,li11tac}.

Because singular systems and hence optimal ensemble controls cannot be derived analytically except for in the simplest cases, an accurate and direct numerical method for approximating the former is a prerequisite for applying this new theory.  Given an appropriate numerical approximation to the singular system $(\sigma_n,\mu_n,\nu_n)$ for the operator $L$ of an ensemble controllable system, the series (\ref{cont1}) can be truncated to synthesize an approximation to $u$ that results in $||X_T-X_F||_K<\eP$ as desired. In addition, a numerical test of the criteria (\ref{contcond}) for ensemble controllability is a natural extension of such a framework, which we present in the following section.

\section{Numerical approximation of ensemble controls} \label{sec_method}

The singular system as defined above is the infinite-dimensional analogue of the well-known singular value decomposition (SVD) for matrices \cite{demmel97}.  A natural approach is therefore to approximate the action of the compact operator $L:\mathcal{H}_T\to\mathcal{H}_K$ in equation (\ref{critdet}) on a function $g\in \mathcal{H}_T$ by a matrix acting on an appropriate vector of sampled values of $g$.  Then the SVD can be used to approximate the singular system of the operator, and thereby also the minimum norm ensemble control $u$.  Let $\{\beta_j\}$ be a finite collection of points that are distributed uniformly throughout the space $K$ and indexed by $j=0,1,2,\ldots,P$, and let $\{t_k\}$ be a collection of points that linearly interpolate the time domain $[0,T]$ for $k=0,1,\ldots,N$, including endpoints, with $t_{k}-t_{k-1}=\delta$.  Using this grid of nodes, we make the approximation
\begin{align} \label{discret}
(Lg)(\beta)&=\int_0^T\Phi(0,t,\beta)B(t,\beta)g(t)dt \nonumber \\
& = \sum_{k=1}^N \int_{t_{k-1}}^{t_k}\Phi(0,t,\beta)B(t,\beta)g(t)dt \nonumber\\
& \approx \sum_{k=1}^N\delta\Phi(0,t_k,\beta)B(t_k,\beta)g(t_k).
\end{align}
The action of the operator $L$ on a function $g\in \mathcal{H}_T$ can be approximated by the action of a block matrix $W\in\mathbb{R}^{nP\times mN}$, with $n\times m$ blocks $W_{jk}=\delta \Phi(0,t_k,\beta_j)B(t_k,\beta_j)$, on a vector $\hat{g}\in \mathbb{R}^{mN}$, with $N$ blocks   $\hat{g}_{k}=g(t_k)$ of dimension $m\times 1$.  If the SVD of this matrix is $W=U\Sigma V^\dagger$, and $\bar{u}_j$ and $\bar{v}_j$ are columns of $U$ and $V$, respectively, corresponding to the singular value $s_j$, then $WW^\dagger \bar{u}_j=s_j^2\bar{u}_j$ and $W^\dagger W\bar{v}_k=s_j^2\bar{v}_k$.  Therefore the SVD $(s_j,\bar{v}_j,\bar{u}_j)$ of the matrix $W$ approximates the singular system $(\sigma_j,\mu_j,\nu_j)$ of the operator $L$, where $\bar{v}_j$ and $\bar{u}_j$ are discretizations of $\mu_j$ and $\nu_j$, respectively.  Now suppose that $\hat{\xi}\in\bR^{nP}$ is given by $\hat{\xi}_k=\xi(t_k)$ for a function $\xi\in\cH_K$.  Then the minimum norm solution $\hat{g}^*$ that satisfies $W\hat{g}=\hat{\xi}$ is given by $\hat{g}^*=W^\dagger z$ where $WW^\dagger z= \hat{\xi}$ \cite{luenberger69}, so applying basic properties of the SVD results in
\begin{align} \label{appsol}
\hat{g}^*=\sum_{j=1}^{mq}\frac{\hat{\xi}^\dagger \bar{u}_j}{s_j}\bar{v}_j.
\end{align}
The components of the synthesized minimum norm control $\hat{u}^*=(\hat{u}_1^*,\ldots,\hat{u}_m^*)^\dagger$ are therefore given by
\begin{align} \label{syncont}
\hat{u}_k^*=\sum_{j=1}^{q}\frac{\hat{\xi}^\dagger \bar{u}_{k+m(j-1)}}{s_{k+m(j-1)}}\bar{v}_{k+m(j-1)}.
\end{align}
Note that the time and parameter discretizations $N$ and $P$ must be chosen such that $nP\leq mN$, so that the pair $(W,\hat{\xi}\,)$ represents an underdetermined system and therefore a minimum norm and not a least squares problem.  The number $q$ of eigenfunctions used in the approximation is limited by $q\leq P$.

The use of a Riemann sum quadrature formula to approximate the action of a Fredholm integral operator of the first type by a matrix, so that the  SVD can be used to approximate the singular system of the operator, has previously been examined as part of a least-squares type method \cite{hanson71}.  An analysis of numerical methods for approximating solutions to Fredholm integrals of the first type, as well as an examination of accuracy and conditioning issues, has also been performed \cite{hansen92}.  The latter work includes a discussion of the Picard criterion,  which asserts that there exists a square integrable solution to the integral equation (\ref{critdet}) only if (i) holds in (\ref{contcond}).  A check of the Picard criterion might provide a test for ensemble controllability, but unfortunately it is very problematic to test this asymptotic property numerically.

The most important computational issue is preventing the aggregation of numerical errors.  These first arise from computation of the flow $\Phi$, which must be done using numerical integration when the system is time-varying.  A practical relative error tolerance for solving ODE systems is $\cO(10^{-6})$.  Another source of numerical error arises from computation of the SVD.  We have found that in order to prevent these errors from dominating the synthesized control, it is appropriate to choose $q$ in (\ref{appsol}) such that the corresponding first and last singular values used satisfy $s_1/s_{mq}<10^4$.

\section{Examples and Discussion} \label{sec_ex}

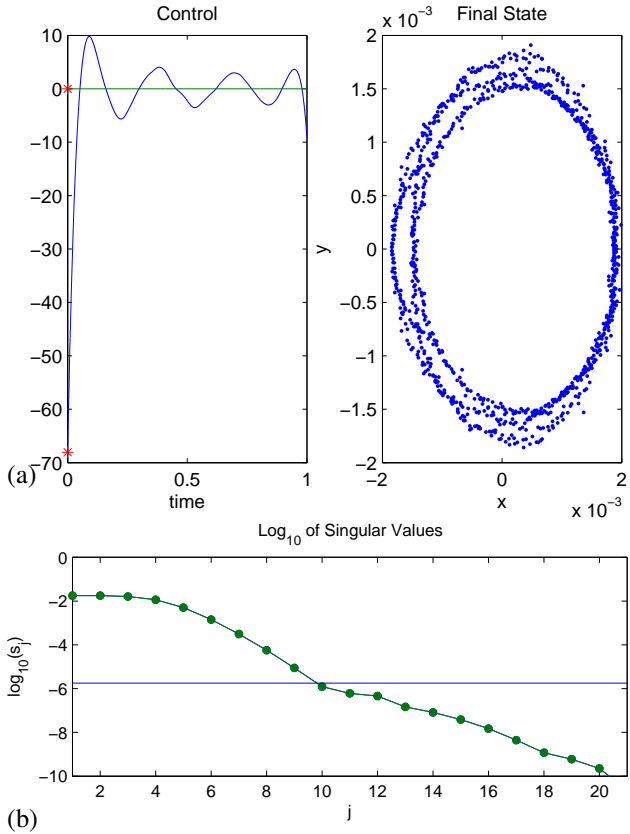
\begin{figure}[t]
\centerline{\includegraphics[width=1.10\linewidth]{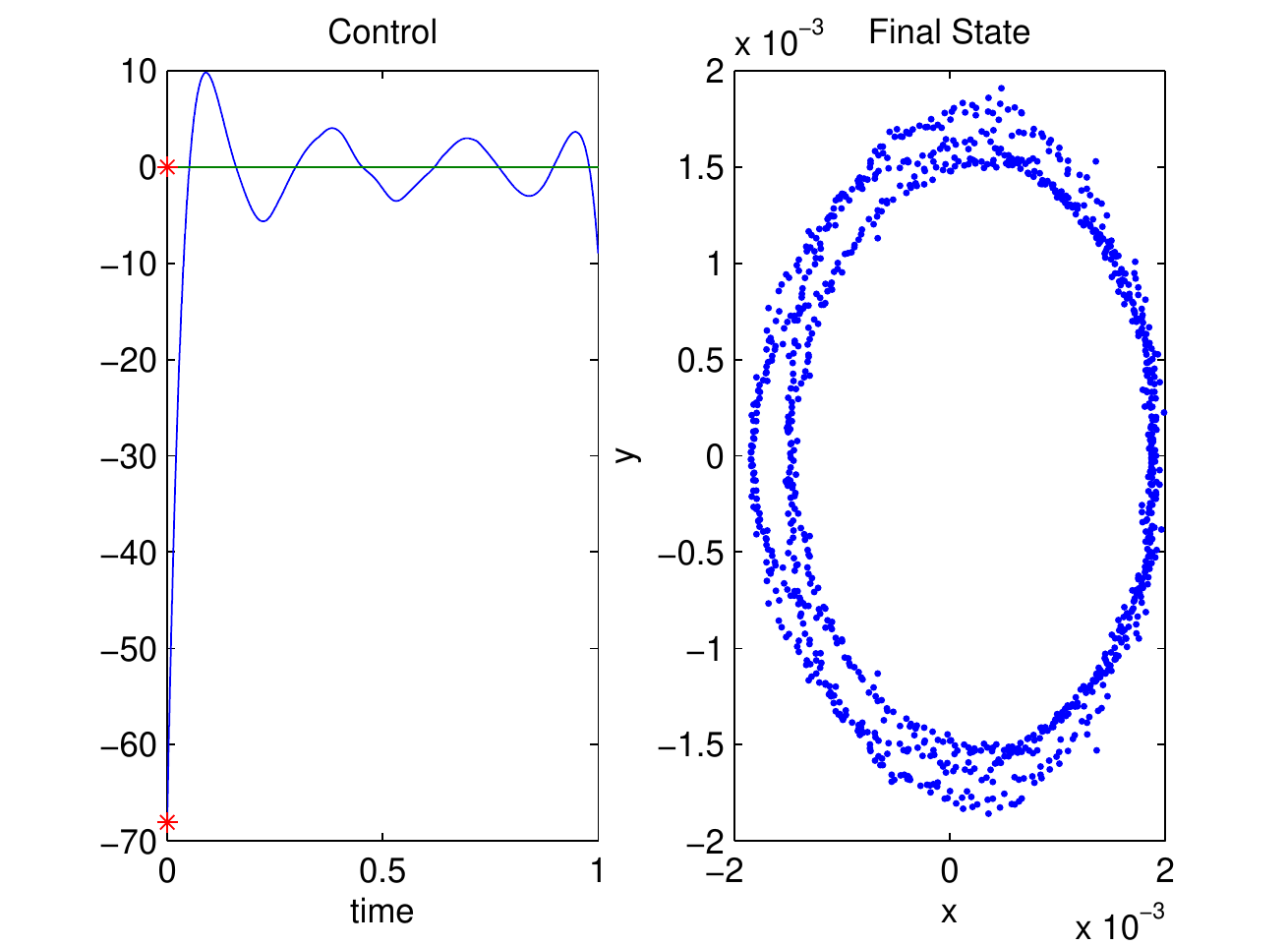}}
\centerline{ \includegraphics[width=1.10\linewidth]{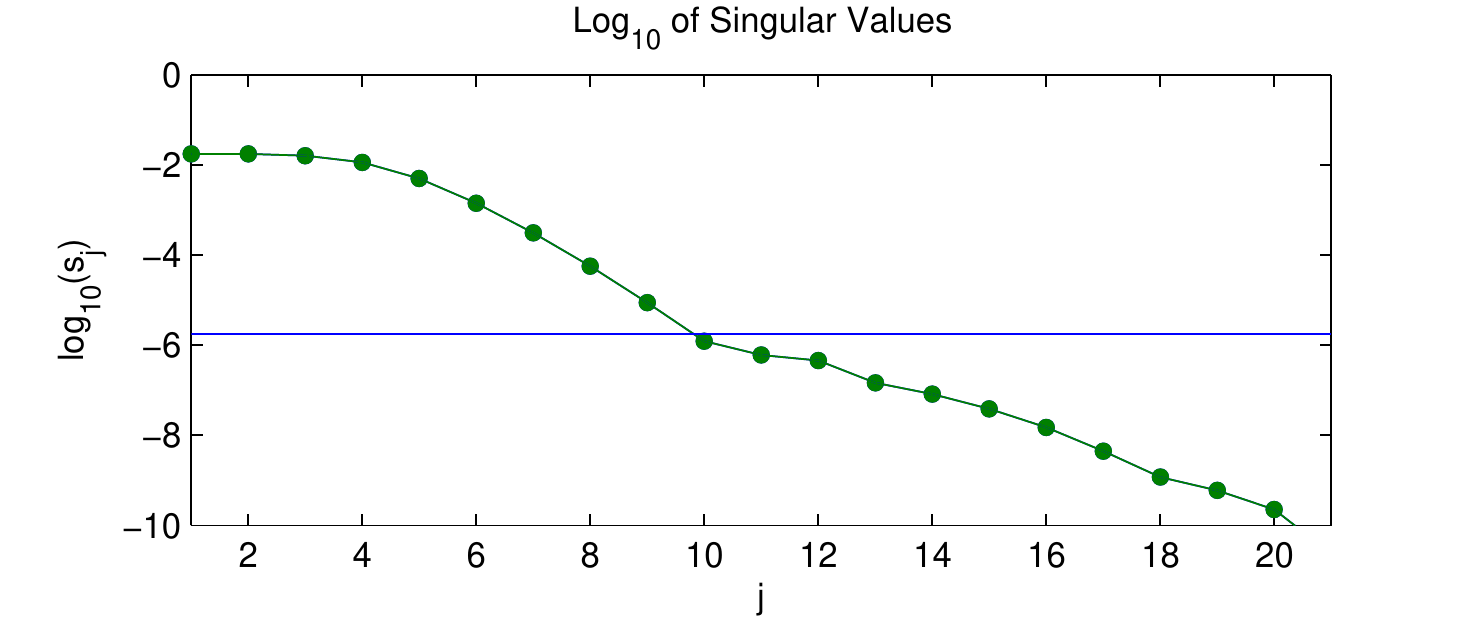}}
\begin{picture}(0,0)
\put(0,140){(a)} \put(0,10){(b)}
\end{picture}
\vskip-.4cm
\caption{Simulation of system ensemble (\ref{harmosc}) for $N=20000$, $T=1$, $P=20$, and $\w\in[-10,10]$.  The initial and target states are $X_0(\w)=(1,0)^\dagger$ and $X_F(\w)=(0,0)^\dagger$.  The matrix $W$ is computed in about 5 seconds, and the SVD is computed in under 1 second. (a) The optimal control law $(u(t),v(t))$ for $t\in[0,1]$ (left), and the final states for all systems $\w\in[-10,10]$.  (b) The singular values $\{s_j\}$ of $W$ on a $\log_{10}$ scale, with the $s_1/s_j<10^4$ cutoff indicated.  Here $q=9$ singular vectors are used to synthesize the control.} \label{figosc}
\end{figure}

In order to illustrate the performance of our method, we will present several examples in this section.  The computations are performed using MATLAB on a PC with a 3.33GHz processor.  We first revisit the optimal control of an ensemble of harmonic oscillators,  which has been previously examined in detail \cite{li11tac}.  This system was proven to be ensemble controllable, and the eigenfunctions of the operator that characterizes the system dynamics are related to the family of prolate spheroidal wave functions.  The dynamics are given by

\begin{figure}[t]
\vskip-.2cm
\centerline{\includegraphics[width=1.05\linewidth]{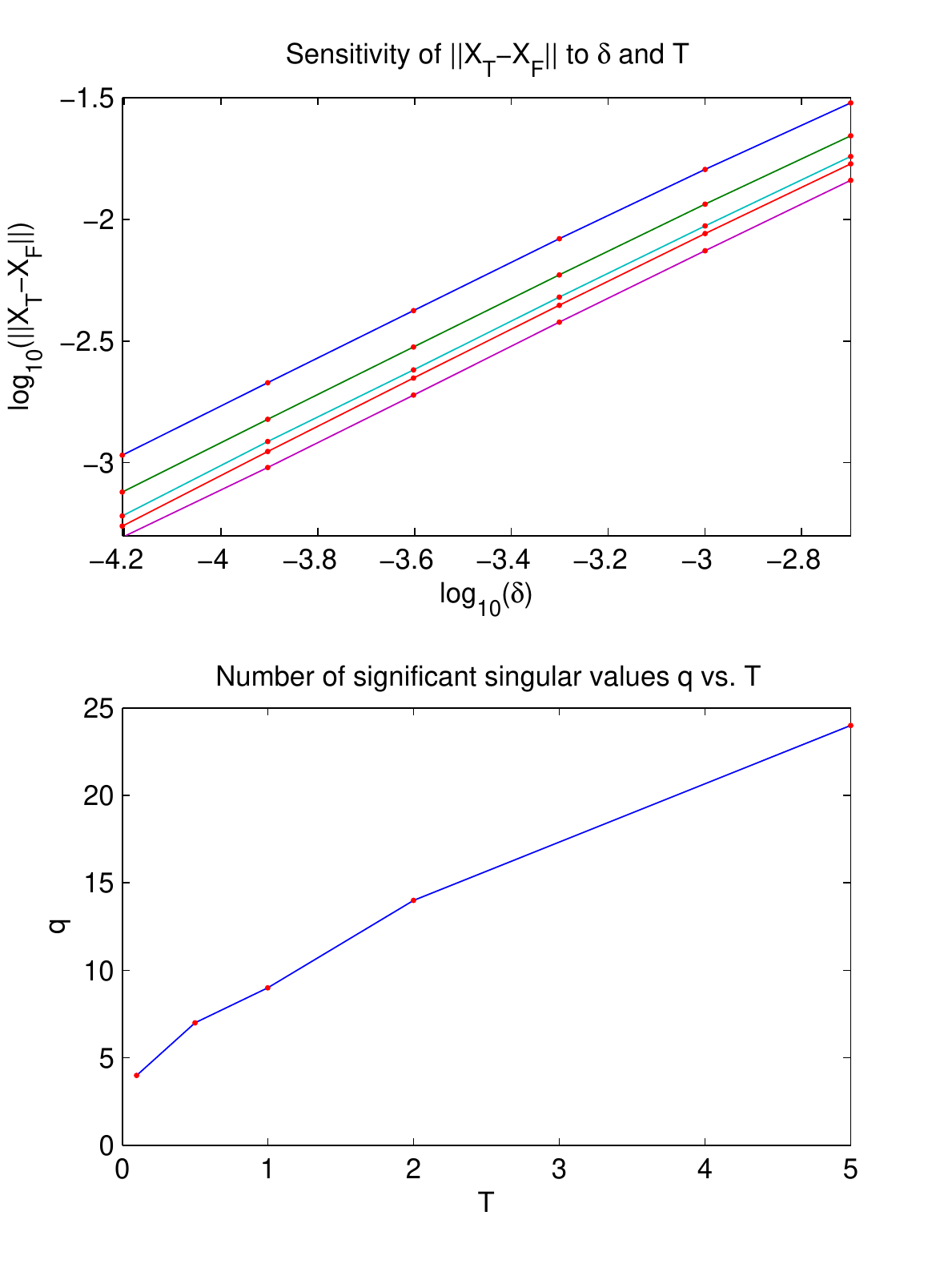}}
\begin{picture}(0,0)
\put(0,195){(a)} \put(0,30){(b)}
\end{picture}
\vskip-.8cm
\caption{The simulation in Figure \ref{figosc} is repeated for different values of time horizon $T$ and time step $\delta$, where $N=T/\delta$, and $P=40$ is used in each instance.  (a) The norm of the error in the final state is plotted as a function of $1/\delta=N/T$. The slope of the lines is very close to 1, so that the error is proportional to $\delta$.  The lines correspond to $T=0.1$, $0.5$, $2$, $1$, and $5$, from top to bottom, hence a longer time horizon does not necessarily result in improvement.  (b) The number of significant singular values $q$ is plotted as a function of $T$.  } \label{fignt}
\end{figure}
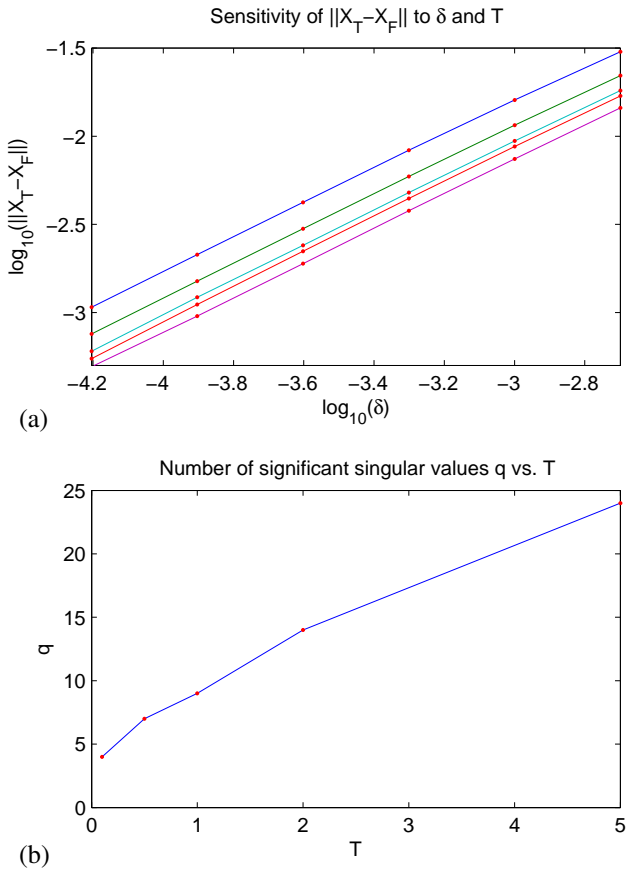
\begin{align} \label{harmosc}
\ddx{t}\bq{\begin{matrix}x(t,\w) \\ y(t,\w)\end{matrix}} = \bq{\begin{matrix} 0 & -\w \\ \w & 0 \end{matrix}} \bq{\begin{matrix} x(t,\w) \\ y(t,\w)  \end{matrix}} + \bq{\begin{matrix} u_1(t) \\ u_2(t)  \end{matrix}},
\end{align}
where $\w\in K=[\w_1,\w_2]\subset\bR$, the instantaneous state is $X(\cdot,\w)=(x(\cdot,\w),y(\cdot,\w))^\dagger\in \cH_K$, and the control vector is $U=(u_1,u_2)^\dagger\in \cH_T$.   We apply the method described in Section \ref{sec_method} to solve an optimal ensemble control problem for the system (\ref{harmosc}), and the results are shown in Figure \ref{figosc}.  A 4th order Runge-Kutta scheme is used to integrate the system using a relative error tolerance of $\cO(10^{-6})$ for adaptive time stepping, and the synthesized control $\hat{u}^*$ is linearly interpolated at the time points requested by the routine.  The control is similar in shape but not equivalent to previous results \cite{li11tac}, and its performance is similar.  In addition, the results of further numerical experiments that determine the sensitivity of the error $\norm{X_T-X_F}_K$ on the choice of time horizon $T$ and time discretization $\delta=T/N$ are shown in Figure \ref{fignt}.  The error in the final state of the ensemble is proportional to the time step $\delta$, which provides a means to calibrate the number $N$ of time discretization points.

\begin{figure}[t]
\centerline{\includegraphics[width=1.05\linewidth]{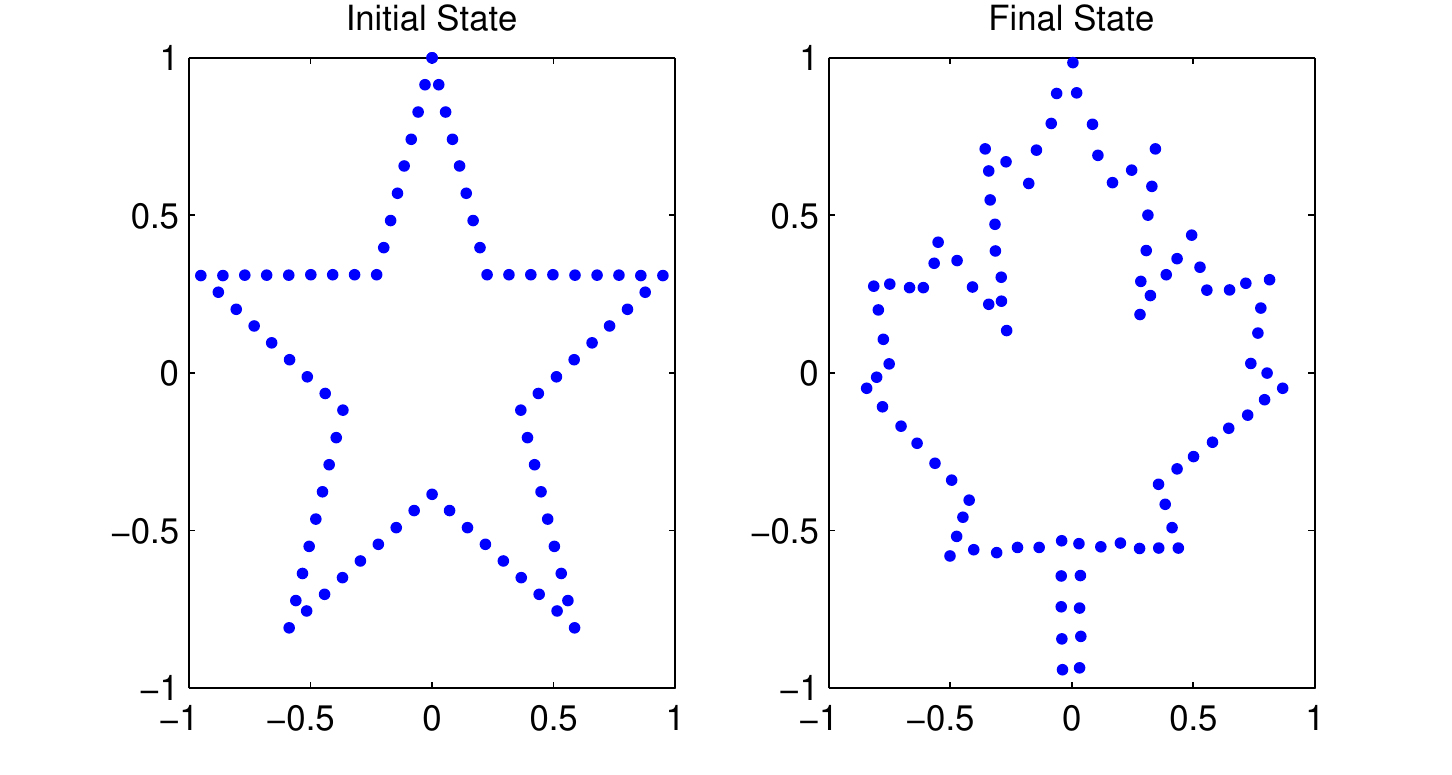}}
\centerline{ \includegraphics[width=1.05\linewidth]{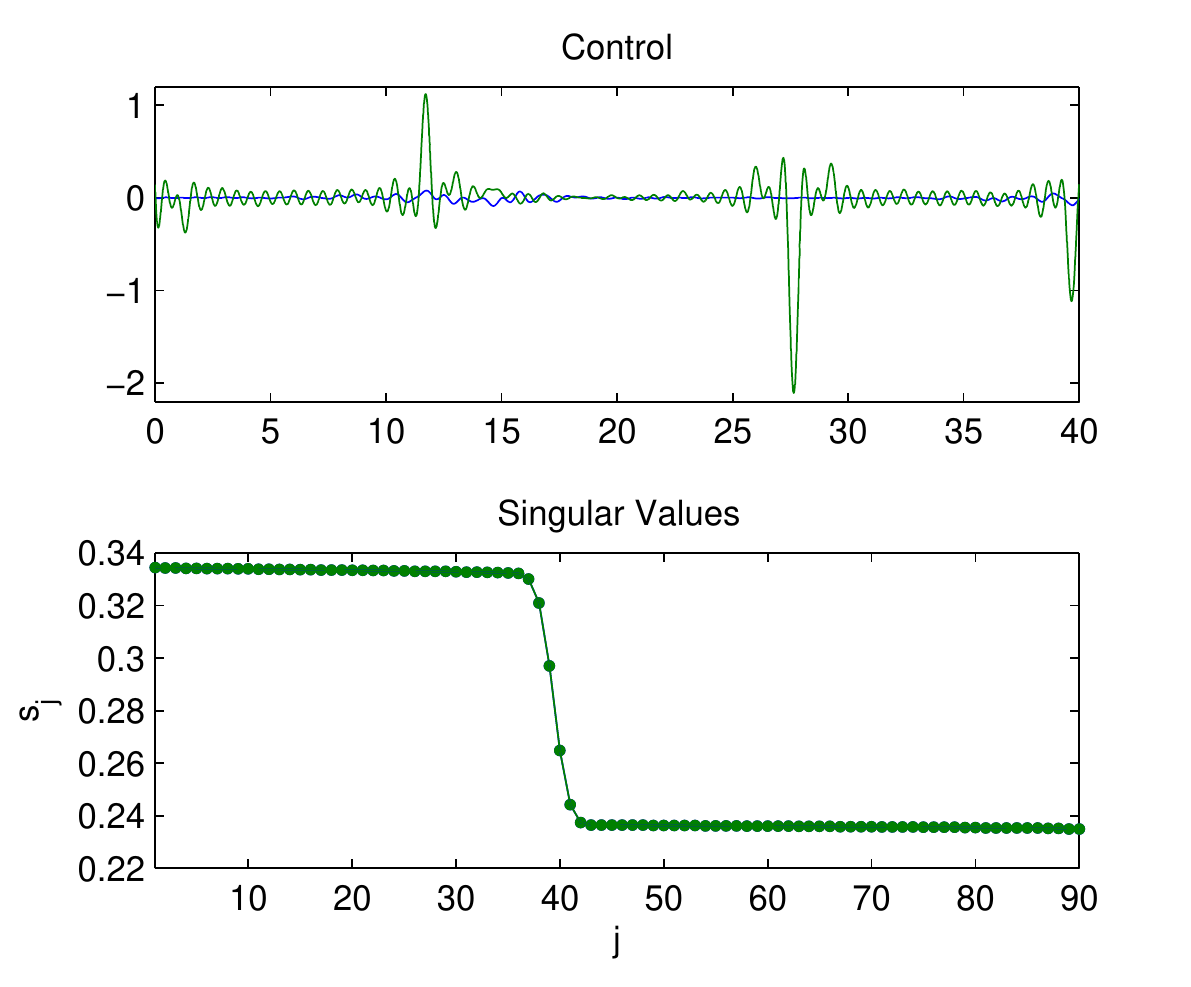}}
\begin{picture}(0,0)
\put(0,230){(a)} \put(0,130){(b)} \put(0,20){(c)}
\end{picture}
\caption{Simulation of system ensemble (\ref{harmosc}) for $N=20000$, $T=40$, $P=89$, and $w\in[-10,10]$, where the initial target states $X_0$ and $X_F$ are arrangements of the $P+1$ oscillators in star and leaf shaped images in the plane, respectively. The average error in the final states of the oscillators is $3.03\times 10^{-4}$, and the maximum error is $0.0167$. (a) Initial and actual final states are plotted. (b) The control that accomplishes the transfer. (c) The spectrum of the SVD is shown. Observe that the spectrum of the SVD differs in form from that which results from the simulation in Figure \ref{figosc}.  Because the conditioning criterion $s_1/s_q<10^4$ is satisfied, all of the singular vectors are used to synthesize the control. } \label{figmaple}
\end{figure}
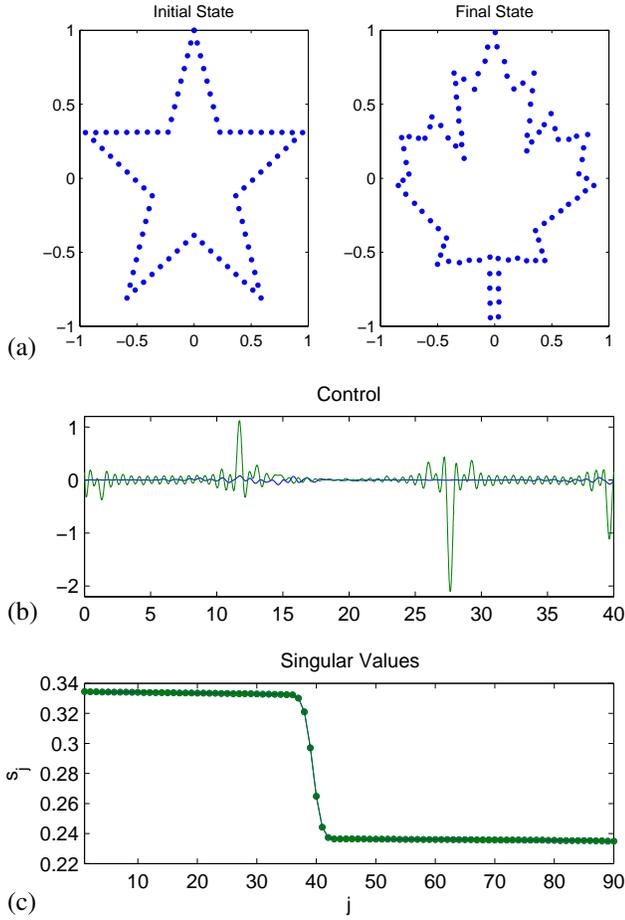

\begin{figure}[t]
\vskip-.7cm
\centerline{\includegraphics[width=1.05\linewidth]{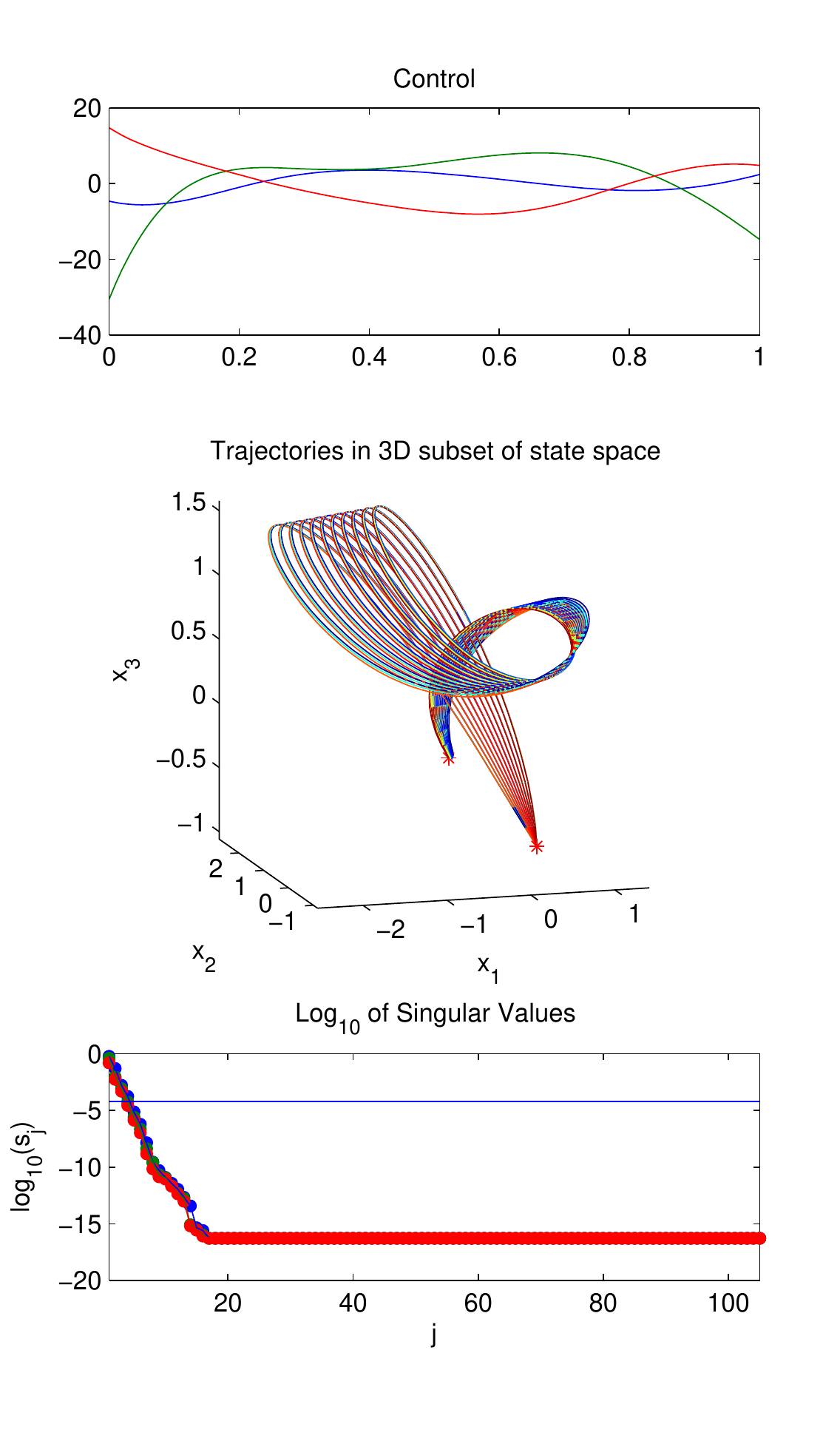}}
\begin{picture}(0,0)
\put(0,330){(a)} \put(0,170){(b)} \put(0,40){(c)}
\end{picture}
\caption{Simulation of the ensemble system (\ref{highsys}) for $N=10000$, $T=1$, $r\in[-.01,.01]$, $c\in[-.1,.1]$, and $P=104$.  The initial and target states are $x_0\approx(0.83,1.38, -1.06, -0.47)^\dagger$ and $x_F=( -0.27, 1.10, -0.28, 0.70)^\dagger$.  The error between the terminal and target states is $||X_T-X_F||_K\approx0.038$.  There are $mq=12$ eigenvectors used to synthesize the control. (a) The minimum norm control. (b) Trajectories in the first 3 coordinates. (c) The spectrum of the SVD decays quickly, as in the example in Figure \ref{figosc}.} \label{fighigh}
\end{figure}
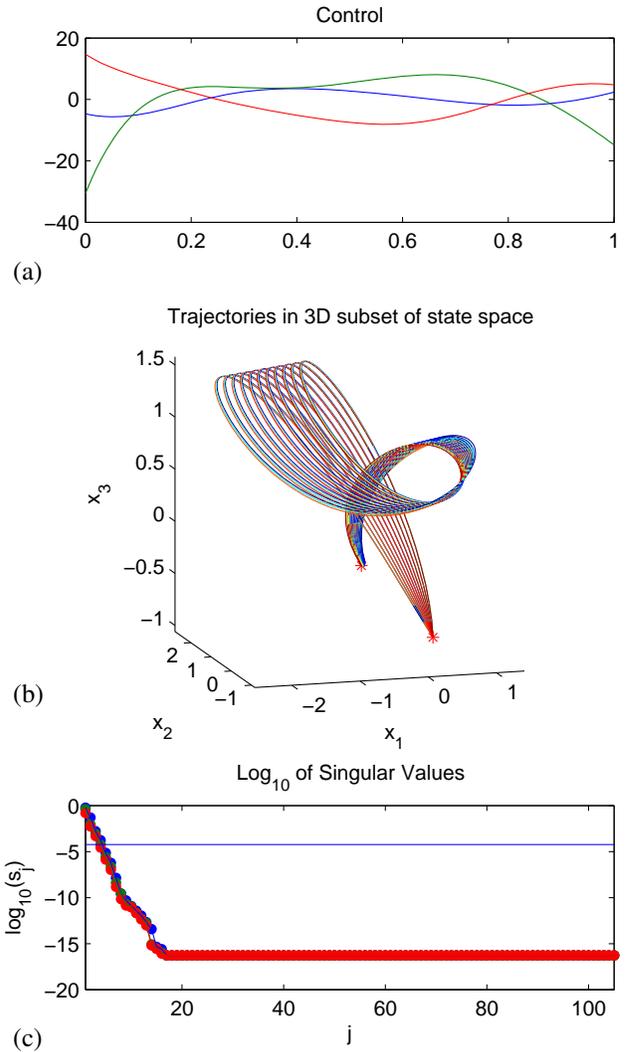

This method can also be used to solve more challenging problems, for example where the initial and target states $X_0$ and $X_F$ are functions of $\w$ in the state space.  We applied our approach to steer an ensemble of harmonic oscillators (\ref{harmosc}) from an initial state arranged in the shape of a star to a target state in the shape of a maple leaf in the plane.  The results of this simulation are shown in Figure \ref{figmaple}, and we encourage the reader to view a video of the transfer that is available online \cite{maplevideo}.  This example is in fact related to a complex problem of importance to the field of NMR in which the initial and target states also depend on system parameters.  In certain experiments specific sub-collections of quantum systems must be excited based on parameter values or the physical position in the sample under study by using so-called selective pulses \cite{kobzar05,pryor07}.  Controls that can create arbitrary patterns in the terminal state of the ensemble as a function of system parameters are therefore desired.

System ensembles governed by dynamics of higher dimension with variation in several parameters can also be controlled using our approach.  Consider for example the time-varying system
\begin{align} \label{highsys}
\dot{X}(t,r,c) & = [A_0+A_1\sin(2\pi t)+A_2 r] X(t,r,c) \nonumber\\
& \quad + [B_0+B_1\tfrac{1}{1+t}+B_2c]u(t),
\end{align}
where $X(t,r,c)\in\bR^4$, $u(t)\in\bR^3$, $r\in[r_1,r_2]$ and $c\in[c_1,c_2]$ are parameters, and $A_0, A_1, A_2 \in\bR^{4\times4}$ and $B_0, B_1, B_2\in\bR^{4\times 3}$ have random entries generated using a standard normal distribution.  Suppose that the initial and target states $X_0(r,c)=x_0\in\bR^4$ and $X_F(r,c)=x_F\in\bR^4$ also similarly randomly generated.  Although ensemble controllability properties of this system are not straightforward to determine, a control synthesis can nevertheless be attempted, with a successful outcome indicating ensemble reachability of the target state from the initial state.  In Figure \ref{fighigh} we provide the results of a simulation where a randomly generated system of the form (\ref{highsys}) is driven between two randomly generated states.  This system is quite sensitive to variation in the parameter $r$, while the algorithm can compensate for dispersion in the parameter $c$ well.

The theoretical treatment in previous work provides a solid foundation for examining ensemble controllability \cite{li11tac}, but a straightforward test for this property is not yet available.  While it is possible to test for ensemble controllability in certain cases by using Lie algebras \cite{li09tac}, the complexity of the systems encountered in many applications makes this problematic in general.  It is important to explore the relationship between ensemble controllability, the Picard criterion, and the singular values of the integral operator (\ref{critdet}) and its matrix analogue $W$.  Investigation in this direction may lead to an implementable numerical test for general ensemble systems.  In addition, a thorough numerical analysis is required to better understand the accuracy and conditioning properties of this approach, in order to predict its performance under various circumstances, and to determine whether a computational test for ensemble controllability is indeed feasible.

\section{Conclusions} \label{sec_conc}

We have introduced an accurate, stable, and computationally efficient numerical method for synthesizing minimum norm ensemble controls for finite-dimensional time-varying linear systems.  By basing our approach on the singular value decomposition (SVD) and discarding all but the most significant singular values, we have guaranteed accuracy and numerical stability of our method, and have leveraged the efficiency of  widely-used numerical routines.  Furthermore, because the SVD is a finite algorithm, our method does not require any additional optimization steps.  We have demonstrated its effectiveness for designing controls for a variety of system ensembles and state transfers under various challenging conditions, including complicated state transfers and high-dimensional time-varying dynamics.  In addition, we have conducted multiple simulations to illustrate the sensitivity of the method with respect to the chosen time step, and the effect of the time horizon on the condition of the problem.

This work forms the basis for a new set of numerical methods for quantum control by providing an approach to designing excitations for guidance of quantum harmonic oscillators.  In our future work, we plan to consider fixed point and contractive properties of integral operator equations to design fast, iterative methods for ensemble control of nonlinear systems.  We will focus in particular on dynamic equations of bilinear form, which govern many widely studied quantum dynamical phenomena.  A computationally efficient numerical scheme to synthesize controls for this type of problem will facilitate the improvement of pulse design for NMR, and would be of immediate use to experimentalists.

\bibliographystyle{unsrt}
\bibliography{acc_ref}

\end{document}